\documentclass[preprint,12pt]{article}

\usepackage{amsfonts}
\usepackage{amsthm}
\usepackage{amsmath}
\usepackage{amsfonts}
\usepackage{latexsym}
\usepackage{amssymb}

\newtheorem{teo}{Theorem}[section]
\newtheorem*{teoo}{Theorem}
\newtheorem{pro}[teo]{Proposition}
\newtheorem{lem}[teo]{Lemma}
\newtheorem{cor}[teo]{Corollary}

\theoremstyle{definition}
\newtheorem{rem}[teo]{Remark}

\def\EOE{\hfill $\blacktriangle$}

\def\bdem{\begin{proof}}
\def\edem{\end{proof}}
\newcommand{\peso}[1]{ \ \ \text{ \rm  #1 } \ \ }

\def\eps{\varepsilon}

\def\la{\lambda}

\def\N{\mathbb{N}}   %Naturales
\def\Z{\mathbb{Z}}   %Enteros
\def\R{\mathbb{R}}   %Reales
\def\C{\mathbb{C}}   %Complejos

\def\cK{\mathcal{K}}

\def\cR{\mathcal{R}}

\def\ese{\mathcal{S}}
\def\ete{\mathcal{T}}
\def\eme{\mathcal{M}}

\def\cU{\mathcal{U}}

\def\h{{\mathcal H}}
\def\b{{\mathcal B}}

\def\ultrar{\cR^\omega}

\newcommand{\pint}[1]{\displaystyle \left \langle\, #1 \, \right\rangle}

\newcommand{\hil}{\mathcal{H}}
\newcommand{\op}{B(\mathcal{H})}
\newcommand{\opc}{\mathcal{B}_0(\mathcal{H})}

\def\fh{ {\mathcal B}_f({\mathcal H})  }
\def\bh{{\mathcal B}({\mathcal H})}

\def\kh{ {\mathcal B}_0(\h) }

\def\u{ {\mathcal U}({\mathcal H}) }

\def\s{\mathcal{S}}

\newcommand{\mat}{\mathcal{M}_n(\mathbb{C})}

\newcommand{\matsa}{\mathcal{H}(n)}
\newcommand{\matu}{\mathcal{U}(n)}

\newcommand{\matinv}{\mathcal{G}\textit{l}\,(n)}

\newcommand{\spec}[1]{\sigma\left(#1\right)}

\newcommand{\sub}[2]{{#1}_{\mbox{\tiny{${#2}$}}}}

\newcommand{\convm}{\xrightarrow[m\rightarrow\infty]{}}

\newcommand{\convsot}{\xrightarrow[n\rightarrow\infty]{\mbox{\tiny{SOT}}}}

\def\puno{\sub{p}{2\pi}^{(n)}}
\def\pdos{\sub{p}{s}^{(n)}}
\def\ptres{\sub{p}{c}}
\def\pcuatro{\sub{p}{0}}

\begin{document}

\title{\huge{\bf{Thompson-type formulae}}}

\author{Jorge Antezana, Gabriel Larotonda and Alejandro Varela}

\date{}

\maketitle

\begin{abstract}
Let $X$ and $Y$ be two $n\times n$ Hermitian matrices. In the article \textit{Proof of a conjectured exponential formula} (Linear and Multilinear Algebra (19) 1986, 187-197)  R. C. Thompson proved that there exist two $n\times n$ unitary matrices $U$ and $V$ such that
$$
e^{i\, X}e^{i\, Y}=e^{i\,  UXU^*+VBV^*}\,.
$$
In this note we consider extensions of this result to compact operators as well as to operators in an embeddable II$_1$ factor.
\footnote{

$\ \,$ MSC[2010]  Primary 47B15,  Secondary 47C15, 47A05
}
\end{abstract}

\tableofcontents

\section{Introduction}

In his 1986 paper \cite{thompson}, studying the product $e^{iX}e^{iY}$ (with $X,Y$ Hermitian matrices) R. C. Thompson considered the analytic map $\xi(w)=e^{iX}e^{iwY}$ defined for some $w\in \mathbb C$ in a neighborhood of the unit interval. Using perturbation theory techniques, he derived a series of inequalities  concerning the eigenvalues of $X,Y$ and those of $Z=log(e^{iX}e^{iY})$. The family of inequalities found by Thompson happened to relate to those proposed by R. Horn \cite{Horn} as the complete solution of the following (seemingly) elementary problem: find necessary and sufficient conditions on the eigenvalues of the Hermitian matrices $A,B,C$ in order to have $UAU^*+VBV^*=C$ for some unitary matrices $U$ and $V$. By that time, V. B. Lidskii had recently published the paper \cite{[L]}, announcing the proof of Horn's conjecture (see Appendix A below for a brief exposition on the subject). Thus, Thompson's computations lead him to conclude that there existed unitary matrices $U,V$ such that $Z=UXU^*+VYV^*$. However, details on Lidskii's proof never saw the light, and for a very long time, Horn's conjecture remained open and consequently, Thompson's result was gently archived. It was not until twelve years later that a proof of Horn's conjecture was given in two exceptional papers, the first one due to A. Klyachko  \cite{[Kl1]} and the second one due to A. Knutson and T. Tao \cite{[KT]}.

Later, Horn's result was extended to the infinite dimensional setting by Bercovici et. al. in two papers  \cite{[BL],[BLT]} that deal with the case of operators in an embeddable II$_1$ factor and with compact operators respectively. Then, it is only natural to ask for extensions of Thompson's formula on adequate infinite dimensional settings. In this paper, we attempt to provide  generalizations to the setting of compact operators, and to the setting of finite von Neumann algebras.  Our motivation stems for the applications of Thompson's identity to the study of the geometry of the Grassmannian manifold when it is endowed with a left-invariant metric induced by a unitarily invariant norm \cite{zhang}.

What follows are the detailed contents of this paper: Section \ref{pre} establishes the background and notation used throughout. In Section \ref{forcompact} we prove that, given $x,y$ compact Hermitian operators on a separable Hilbert space, there exist unitary operators $u,v$ and an isometry $w$ such that
$$
e^{iwxw^*}e^{iwyw^*}=e^{iuwxw^*u^*+ivwyw^*v^*}.
$$
This is the content of Theorem \ref{porron}. The strategy adopted for our proof is the following: cut $x,y,z$ with adequate finite rank projections, apply Thompson's formula to produce  sequences $(u_nxu_n^*)$, ($v_nyv_n^*$). Then, the main task is to extract convergent subsequences of the exponents obtained. 
 The isometry $w$ of the above formula is there to compensate for the possible variations on the kernels of the operators involved: recall one of the consequences of a result by D. Voiculescu \cite{voicu}, that tells us that the unitary orbit $\{uxu^*\}$ of a bounded operator $x$ is closed if and only if $x$ is a finite rank operator. It is however unclear if the isometry $w$ introduced is an artifact due to the strategy of proof we adopted, or there is something deeper going on underneath.

Regarding the final sections of the paper, in Section \ref{fini} we prove that in an adequate von Neumann algebra (i.e., in an embeddable II$_1$ factor), it is possible to obtain an approximate, but uniform, Thompson formula: the precise statement is given in Theorem \ref{Thompson en el ultranoseque}. Here, we use the existence of a convenient matrix approximation, Thompson's formula in finite dimensions, and the extension of Horn's result to the setting of embeddable II$_1$ factors. We finish with an Appendix that contains a brief survey of Horn's conjecture, its solution, and recent generalizations to the settings considered here.

\section{Preliminaries}\label{pre}

Let ${\mathcal H}$ be a complex separable and infinite dimensional Hilbert space. In this paper $\bh$, $\opc$ and $\fh$ stand for the sets of bounded linear operators, compact operators and finite rank operators in ${\mathcal H}$ respectively. The unitary group of $\bh$ is indicated by $\u$. If $x\in \bh$, then $\|x\|$ stands for the usual uniform norm, and we will use $|\cdot|$ to indicate the modulus of an operator, i.e. $|x|=\sqrt{x^*x}$. We indicate with  $\bh_{h}$ (resp. $\kh_h$) the real linear space of Hermitian elements (resp. Hermitian compact elements) of $\bh$. Given $\eta,\zeta\in\hil$, by means of $\eta\otimes\zeta$ we denote the rank one operator defined by $\eta\otimes\zeta(\xi)=\pint{\xi,\zeta}\eta$.

\medskip

On the other hand, throughout this paper $\mat$ denotes the algebra of complex $n\times n$ matrices, $\matinv$ the group of all invertible elements of $\mat$,  $\matu$ the group
of unitary $n\times n$ matrices, and $\matsa$ the real subalgebra of Hermitian matrices.

\medskip

Given $x\in\bh$ (or $T \in \mat$), $R(x)$  the range or image of $x$, $N(x)$ the null space of $x$, and $\spec{x}$ denotes the spectrum of $x$. If $x$ is normal (i.e. $xx^*=x^*x$), then $E_x(\Omega)$ denotes the spectral measure of $x$ associated to the (measurable) subset $\Omega$ of the complex plane. 

\medskip

Let us give a precise statement of Thompson's formula:

\begin{teoo}[Thompson]\label{thompson}
Given $X,Y\in\matsa$, there exist unitary matrices $U,V\in\mat$ such that
$$
e^{i\, X}e^{i\,Y}=e^{i (UXU^*+VBV^*)}\,.
$$
\end{teoo}

\subsection{Some preliminaries on II$_1$ factors}

  Throughout this section, $\ultrar$ denotes the ultrapower the hyperfinite II$_1$ factor, and $\eme$ denotes any II$_1$ factor that can be embedded in $\ultrar$. We are going to use the Greek letter $\tau$ to denote the normalized tracial state of $\eme$. Given an Hermitian element $a\in\eme$, it can be written as
$$
a=\int_0^1 \la_a(t)\ de(t)\,,
$$
where $\la_a$ is a non-increasing right continuous function, and $e(\cdot)$ is a spectral measure on $[0,1)$ such that $\tau(e(t))=t$.

\medskip

  One of the characterizations of embeddable factors is the existence of a ``sequence of matricial approximations'' for any finite family of Hermitian elements. This notion is described more precisely in the following theorem:

\begin{teo}\label{aprox}
Let $a_1,\ldots,a_k$ be Hermitian elements of $\ultrar$. Then, there are integer numbers $1\leq n_1<n_2<\ldots$ and Hermitian matrices
$X^{(m)}_1,\ldots,X^{(m)}_k\in\eme_{n_m}(\C)$ such that for every non-commutative polynomial $p$ it holds that
$$
\tau\big(p(a_1,\ldots,a_k)\big)=\lim_{m\to\infty} \tau_{n_m} \big(p(X^{(m)}_1,\ldots,X^{(m)}_k)\big)
$$
where $\tau_{n_m}$ is the normalized trace of $\eme_{n_m}(\C)$. Moreover, the matrices can be taken so that for each $j\in\{1,\ldots k\}$ we have that $\|X^{(m)}_j\|\leq \|a_j\|$ for every $m\in\N$.
\end{teo}

  Given a matrix $M\in\mat$ whose eigenvalues arranged in non-increasing order are denoted by $\la_1,\ldots,\la_n$, let $\la_M$  denote the real valued function defined in $[0,1)$ in the following way:
$$
\la_M(t)=\sum_{j=1}^n \la_j \chi_{[\frac{j-1}{n},\frac{j}{n})}\,.
$$
With this notation, the following result is a direct consequence of Theorem \ref{aprox}, and the reader is referred to \cite{[BL]} for a detailed proof:

\begin{cor}
Let $a\in\eme$ be an Hermitian element, and $\{X^{(m)}\}_{m\in\N}$ a sequence of matricial approximations of $a$. Then, $\la_{X^{(m)}}\convm \la_a$ almost everywhere.
\end{cor}

  Finally, we mention the following result valid in every finite factor:

\begin{pro}[Kamei \cite{Kamei}]\label{Kamei}
Let $a$ and $b$ be Hermitian elements of a finite factor $\eme$. Then, the following statements are equivalent:
\begin{enumerate}
	\item $\la_a=\la_b$;
	\item $a$ belongs to the norm closure of the unitary orbit of $b$.
\end{enumerate}
\end{pro}

\section{Thompson-type formulae for compact operators}\label{forcompact}

  Throughout this section, given a compact operator $x$, the eigenvalues of $x$ are arranged in non-decreasing order with respect to their moduli, i.e., if $i\leq j$ then $|\la_i(x)|\geq |\la_j(x)|$. 

\begin{teo}\label{triad}
Given $x, y\in \kh_h$, there exist unitary operators $u_k$ and $v_k$ $\in\b(\h)$, for $k\in\mathbb N$, such that
\begin{equation}\label{eq triad}
e^{i\, x}e^{i\, y}=\lim_{k\to\infty}e^{i\,  u_kxu_k^*+i\, v_kyv_k^*}
\end{equation}
\end{teo}

\bdem
Let $x=u|x|$ and $y=v|y|$ be polar decompositions of $x$ and $y$, and
$$
|x|=\sum_{j\in\N} \la_j(|x|) \beta_j \otimes \beta_j\,\peso{and} |x|=\sum_{j\in\N} \la_j(|y|) \zeta_j \otimes \zeta_j\,.
$$
spectral decompositions of $|x|$ and $|y|$ respectively. Recall that the eigenvalues are arranged in non-increasing order. Define
$$
x_k=\sum_{j=1}^k \la_j(|x|) \beta_j \otimes (u\beta_j)\,\peso{and} y_k=\sum_{j=1}^k \la_j(|y|) \zeta_j \otimes (v\zeta_j)\,,
$$
and $\s_k=R(x_k)+ R(y_k)$. Then $x_k(\s_k)\subset \s_k$ and $y_k(\s_k)\subset \s_k$. So, $x_k, y_k\in\b(\s_k)\simeq M_n(\mathbb C)$ (where $n=dim(\s_k)$). On the other hand,
\begin{equation}\label{conv unif}
e^{ix}e^{iy}=\lim_{k\to\infty} e^{ix_k}e^{iy_k}\,.
\end{equation}

Due to Thomspon's formula for matrices, there exist $u_k, v_k$ unitary linear transformations in $\s_k$ (which means that $u_k u_k^*=p_{\s_k}$ and $v_k v_k^*=p_{\s_k}$, where $p_{\s_k}$ denotes the orthogonal projection onto $\s_k$) such that
\begin{equation}\label{formula thompson en sk}
e^{i\,x_k}e^{i\,y_k}= e^{i\, u_k x_k u_k^*+i\, v_k y_k v_k^*}.
\end{equation}

  We can extend $u_k, v_k\in\b(\s_k)$ to the unitaries $\tilde u_k=u_k+p_{\s_k^\perp}$ and $\tilde v_k=v_k+p_{\s_k^\perp}$ $\in\b(\h)$. Then from the equality (\ref{formula thompson en sk}) valid in $\s_k$ we get the following in $\b(\h)$

\begin{equation}\label{formula thompson en B(H)}
e^{i\,x_k}e^{i\,y_k}= e^{i\, \tilde u_k x_k \tilde u_k^*+i\, \tilde v_k y_k \tilde v_k^*}.
\end{equation}
Since $(\tilde u_k x\tilde u_k^*+\tilde v_k y\tilde v_k^*)- (\tilde u_k x_k\tilde u_k^*+\tilde v_k y_k\tilde v_k^*)\to 0$, using \eqref{conv unif} we get
$$
e^{i\, (\tilde u_k x\tilde u_k^*+\tilde v_k y\tilde v_k^*)}- e^{i\, x}e^{i\, y}\xrightarrow[k\rightarrow\infty]{\|\,\cdot\,\|} 0.
$$
\edem

%Given an unitary operator $u$, there always exists an Hermitian operator $x$ such that $e^{ix}$ and $\|x\|\leq \pi$. So, an important particular case of Theorem \ref{triad}, is when $\|x\|\leq \pi$ and $\|y\|\leq \pi$. Under these additional hypotheses, the next result can be obtained.

Since the unitary orbit of a fixed operator in $\bh$ is not closed in general, to avoid the limit in \eqref{eq triad} we have to pay some price. The following theorem follows this path.

\begin{teo}\label{porron}
Given $x, y\in \kh_h$, there is an isometry $w\in\op$, and unitary operators $u$ and $v$ such that
$$
e^{i\, wxw^*}e^{i\, wyw^*}=e^{i\,  u(wxw^*)u^*+i\, v(wyw^*)v^*}\,.
$$
\end{teo}

%\begin{rem}
%The condition $z\in\kh_h$ can be deduced from the other hypothesis. Indeed, if $x,y\in\kh_h$, then $e^{iz}-1=e^{ix}e^{iy}-1$ is compact. Then, as $\|z\|\leq \pi$, the operator $z^{-1}(e^{iz}-1)$ is bounded and invertible. As
%$$
%z=(e^{iz}-1)\Big(z^{-1}(e^{iz}-1)\Big)^{-1},
%$$
%we deduce that $z$ have to be also compact.
%\end{rem}

\begin{rem}\label{version equivalente} Another way to state the theorem follows: there is a bigger Hilbert space $\cK$ containing $\hil$ such that the extensions $\widehat{x}, \widehat{y}\in B(\cK)$  defined by
$$
\widehat{x}=\begin{pmatrix}
x&0\\
0&0
\end{pmatrix}
\begin{array}{l}
\hil\\
\cK\ominus\hil
\end{array}\ ,\ \
\widehat{y}=\begin{pmatrix}
y&0\\
0&0
\end{pmatrix}
\begin{array}{l}
\hil\\
\cK\ominus\hil
\end{array}
$$
satisfy the identity
$
e^{i\, \widehat{x}}e^{i\, \widehat{y}}=e^{i\, (u\widehat{x}u^*+ v\widehat{y}v^*)}\,,
$
for some unitary operators $u$ and $v$ acting on $\cK$. \EOE
\end{rem}

  Let us roughly sketch the idea behind the proof. We know that there are unitary operators $u_n,v_n\in\cU(\hil)$ such that
$$
e^{ix}e^{iy}=\lim_{n\to\infty}e^{i(u_nxu_n^*+v_nyv_n^*)}\,.
$$
Let $z_n=u_nxu_n^*+v_nyv_n^*$. Extending to a  bigger space $\cK$ the operators $z_n$, $u_n$, $v_n$, $x$ and $y$ as in the previous remark, we can conjugate the sequence $\{\widehat{z}_n\}_{n\in\N}$ with unitary operators $w_n$ acting on $\cK$ so that $e^{\widehat{z}_n}=e^{w_n\widehat{z}_nw_n^*}$, and the modified sequence $\{w_n\widehat{z}_nw_n^*\}_{n\in\N}$ has a convergent subsequence.  If $\widehat{s}$ denotes the limit of that subsequence, provided  $\dim\cK\ominus\hil=\infty$, we can always find two unitary operators $\widehat{u}_0$ and $\widehat{v}_0$  such that
$$
\widehat{s}=\widehat{u}_0\,\widehat{x}\,\widehat{u}_0^*+\widehat{v}_0\,\widehat{y}\,\widehat{v}_0^*.
$$
As this limit $\widehat{s}$ satisfies that $e^{i\widehat{x}}e^{i\widehat{y}}=e^{i\widehat{s}}$, this would complete the proof. Since the proof of Theorem \ref{porron} is rather long, some technical parts are included in the next three lemmas:

%The technical obstacle that main problem is that we have not been able to prove that $\{z_n\}_{n\in\N}$ has a convergent subsequence. That is the reason why we have to extend the operators $x$, $y$, and $z$ to a bigger Hilbert space, where we can solve the technical problems appearing in $\hil$.

\begin{lem}\label{papo}
Let  $\{a_n\}_{n\in\N}$ be a bounded sequence of finite rank normal operators, and let $p_n$ denote the orthogonal projection onto $R(a_n)$. If there exists a finite rank projection $p$ such that $p_n\xrightarrow[n\to\infty]{\|\cdot\|} p$, then $\{a_n\}_{n\in\N}$ has a convergent subsequence.
\end{lem}
\bdem
Since $p_n\xrightarrow[n\to\infty]{\|\cdot\|} p$, the operators $s_n:=p_np+(1-p_n)(1-p)$ converge to $1$ as $n\to\infty$. We can suppose that for every $n\in\N$, $s_n$ is invertible. Note also that $p_ns_n=s_np$. For each $n\in\N$, let $s_n=u_n|s_n|$ be the polar decomposition of $s_n$. Then, straightforward computations show that $p_nu_n=u_np$. So, as the sequence $\{a_n\}_{n\in\N}$ is bounded, $\{u_n^*a_nu_n\}_{n\in\N}$ is a bounded sequence of normal operators whose range is the finite dimensional subspace $R(p)$. Therefore, it has a norm-convergent subsequence. Since $u_n \xrightarrow[n\to\infty]{\|\cdot\|} 1$, the original sequence $\{a_n\}_{n\in\N}$ also has a convergent subsequence.
\edem

\begin{lem}\label{rodriguez}
Let $z\in \kh_h$ be such that $\|z\|\leq\pi$, and let $\{w_n\}_{n\in\N}$ be a bounded sequence of Hermitian compact operators which satisfies:
\begin{enumerate}
	\item[a.)] $e^{iw_n}\xrightarrow[n\to\infty]{\|\cdot\|} e^{iz}$;
	%\item[b.)] $\displaystyle \sup_{n\in \N}\|w_n\|< 2\pi$. 
	\item[b.)] There exists $n_0\in\mathbb N$ and $\eps>0$ such that
	$$\Big(\bigcup_{n\ge n_0} \sigma(w_n) \Big)\cap \Big(\bigcup_{k\in\Z,k\neq 0} (2k\pi-\eps,2k\pi+\eps)\Big)=\varnothing$$
\end{enumerate}
Then, $\{w_n\}_{n\in\N}$ has a convergent subsequence.
\end{lem}

\bdem
Since $\|z\|\leq\pi$ and the operators $w_n$ satisfy condition (b), there exists $\eps>0$ such that it is not contained neither in the spectrum of any $w_n$ nor in the spectrum of $z$, and it satisfies
\begin{eqnarray}
p_n=&E_{e^{w_n}}\Big(B_1\Big(2\sin\frac{\eps}{2}\Big)\Big)&=E_{w_n}((-\eps,\eps))\nonumber\\
p=&E_{e^{z}}\Big(B_1\Big(2\sin\frac{\eps}{2}\Big)\Big)&=E_{z}((-\eps,\eps))\,, \nonumber
\end{eqnarray}
where $B_\alpha(\rho)$ denotes the ball in $\C$ of radius $\rho$ centered at $\alpha$. Standard arguments of functional calculus imply that $p_n\xrightarrow[n\to\infty]{\|\cdot\|} p$. If $\log$ denotes the principal branch of the complex logarithm, then
$$
\log\big((1-p_n)+p_ne^{w_n}\big)=p_nw_n \peso{and} \log\big((1-p)+pe^{z}\big)=pz.
$$
So, $p_nw_n\xrightarrow[n\to\infty]{\|\cdot\|} pz$ because the sequence $\{(1-p_n)+p_ne^{w_n}\}_{n\in\N}$ converge in the norm topology to 
$(1-p)+pe^{z}$, and the holomorphic functional calculus is continuous with respect to this topology. On the other hand, if $q_n=1-p_n$, the sequence $\{w_nq_n\}_{n\in\N}$ satisfies the conditions of Lemma \ref{papo}. Hence, it has a convergent subsequence $\{w_{n_k}q_{n_k}\}_{k\in\N}$. Therefore, $\{w_{n_k}\}_{k\in\N}$ converges, which concludes the proof.
\edem

  The next lemma is a variation of Lemma 4.3 in \cite{[BLT]}, and its proof follows essentially in the same lines. We include a sketch of its proof for the sake of completeness.

\begin{lem}\label{decadentes}
Let $x,y\in\opc_h$, and suppose there exist unitary operators $u_k$ and $v_k$, for $k\in\N$ such that
$$
s=\lim_{k\to\infty} u_k xu_k^*+v_k yv_k^*\,,
$$
for some $s\in\opc$. Then, there exist compact operators $\bar{s}$, $\bar{x}$, $\bar{y}$  satisfying $\bar{s}=\bar{x}+\bar{y}$,  $\spec{\bar{s}}=\spec{s}$, $\spec{\bar{x}}=\spec{x}$, and $\spec{\bar{y}}=\spec{y}$ with the same multiplicity for every non-zero eigenvalue.
\end{lem}
\bdem[Sketch of proof.]  Let $x_k=u_k xu_k^*$, $y_k=v_k yv_k^*$, and $s_k=x_k+y_k$. For each $k\in\N$ consider an increasing sequence of projections $\{p_{k,n}\}_{n\in\N}$ such that $\dim R(p_{k,n})=n$, $p_{k,n}\convsot 1$, and
$$
\eps_n:=\sup_{k\in\N} \big(\|(1-p_{k,n})s_k\|+\|(1-p_{k,n})x_k\|+\|(1-p_{k,n})y_k\|\big) \xrightarrow[n\to\infty]{}0\,.
$$
This last requirement can be achieved by choosing the projections in such a way that they capture for each $n$ as many eigenvectors of $x_k$ and $y_k$ as it is possible, among those corresponding to the biggest eigenvalues (in modulus) of $x_k$ and $y_k$.

Now, consider a fixed increasing sequence of projections $\{q_n\}_{n\in\N}$ such that $\dim R(q_n)=n$, and $q_{n}\convsot 1$, and for each $k\in\N$ define a unitary operator $w_k$ such that
$$
w_k p_{k,n} w_k^*=q_n\,.
$$
Let $\bar{s}_k=w_ks_kw_k^*$, $\bar{x}_k=w_kx_kw_k^*$, and $\bar{y}_k=w_ky_kw_k^*$. Straightforward computations show that these operators satisfy the following inequalities:
\begin{align}
\|\bar{s}_k-q_n\bar{s}_kq_n\|\leq 2\eps_n\,, \,
\|\bar{x}_k-q_n\bar{x}_kq_n\|\leq 2\eps_n\,, \,  \mbox{and} \
\|\bar{y}_k-q_n\bar{y}_kq_n\|\leq 2\eps_n\,. \label{ufff}
\end{align}
Note that, for each $n\in\N$, set $\{q_n\bar{s}_kq_n:\ k\in\N\}$ is bounded, hence totally bounded. So, the first inequality of \eqref{ufff} implies that the set $\{\bar{s}_k:\ k\in\N\}$ is totally bounded as well. Therefore, passing to a subsequence if necessary, we may assume that  the sequence $\{\bar{s}_{k}\}$ converges to a compact Hermitian operator $\bar{s}$. The same argument can be applied to the sequences $\{\bar{x}_k\}$ and $\{\bar{y}_k\}$, and
we get the operators $\bar{x}$, and $\bar{y}$, respectively. Clearly these operators satisfy
$$
\bar{s}=\bar{x}+\bar{y}\,,
$$
and standard arguments of functional calculus show that $\spec{\bar{s}}=\spec{s}$, $\spec{\bar{x}}=\spec{x}$, and $\spec{\bar{y}}=\spec{y}$ with the same multiplicity for every non-zero eigenvalue.
\edem

\bdem[Proof of Theorem \ref{porron}]
Let $z$ be any bounded and Hermitian operator such that $e^{iz}=e^{ix}e^{iy}$. For simplicity, we are going to prove the alternative version of the statement described in Remark \ref{version equivalente}, and without lost of generality, we are going to assume that $\|z\|\leq\pi$. Then note that, since $e^{iz}-1=e^{ix}e^{iy}-1$ and the right hand is compact, then an elementary argument using the funcional calculus of the entire map $F(\lambda)=(e^{i\lambda}-1)\lambda^{-1}$ shows that $z$ is also a compact operator.

By Theorem  \ref{triad}, there are unitary operators $u_n$ and $v_n$ such that:
$$
e^{iz}=\lim_{n\to\infty} e^{i(u_nxu_n^*+v_nyv_n^*)}\,.
$$
Let $z_n:=u_nxu_n^*+v_nyv_n^*$. Since $x$ and $y$ are compact, there exists $M>0$ big enough such that for every $j\geq M$ and every $n\in\N$ it holds that $\la_j(|z_n|)<\pi$. For technical reasons, passing to a subsequence if necessary, we can assume that $\{\la_j(z_n)\}_{n\in\N}$
converges for every $j\in\{1,\ldots,M\}$. Define 
\begin{align*}
\Omega&=\{m\in\N:\ \limsup_{n\to\infty} \la_m(|z_n|)=2k\pi \ \mbox{for some $k\in\N$}\}\\
&=\{m\in\N:\ \lim_{n\to\infty} \la_m(|z_n|)=2k\pi \ \mbox{for some $k\in\N$}\}.
\end{align*}
The second equality holds because $\#\Omega<M$. Let $\{\zeta_j^{(n)}\}_{j\in\N}$ be an orthonormal basis of $\hil$ such that $\zeta_j^{(n)}$ is an eigenvector of $\la_j(z_n)$. Then
\begin{equation}\label{convergen a 2pi}
\lim_{n\to\infty} \pint{|z_n|\zeta_j^{(n)},\zeta_j^{(n)}}=2k\pi\quad \quad\,\ \mbox{for  $j\in\Omega$, and some $k\in\Z$.}
\end{equation}
Let $\cK=\hil\oplus \hil$,  and extend $x$,$y$, $z$ to $\cK$ as:
$$
\widehat{x}=\begin{pmatrix}
x&0\\
0&0
\end{pmatrix}
\begin{array}{l}
\hil\\
\hil
\end{array}\ ,\quad
\widehat{y}=\begin{pmatrix}
y&0\\
0&0
\end{pmatrix}
\begin{array}{l}
\hil\\
\hil
\end{array}\ ,\peso{and}
\widehat{z}=\begin{pmatrix}
z&0\\
0&0
\end{pmatrix}
\begin{array}{l}
\hil\\
\hil
\end{array}\,.
$$
The unitary operators $u_n$ and $v_n$ are also extended, but in this case as the identity in the second copy of $\hil$. Denote with $\widehat{u}_n$ and $\widehat{v}_n$ these extensions. With these definitions, we get
$$
\widehat{z}_n=\widehat{u}_n\,\widehat{x}\,\widehat{u}_n^*+\widehat{v}_n\,\widehat{y}\,\widehat{v}_n^*=\begin{pmatrix}
z_n&0\\
0&0
\end{pmatrix}
\begin{array}{l}
\hil\\
\hil
\end{array}\,.
$$
Fix an orthonormal basis $\{\beta_j\}_{j\in\N}$  of $\hil$, and define for each $n\in\N$ the unitary operator $w_n$ as the unique unitary operator in $B(\cK)$ that satisfies
\begin{align*}
w_n(\zeta_j^{(n)}\oplus 0)&=\ 0\oplus\beta_j\quad\  \mbox{if $j\in\Omega$}\,,\\
w_n(0\oplus\beta_j\ )  &=\zeta_j^{(n)}\oplus 0 \quad \mbox{if $j\in\Omega$}\,,\\
w_n(\zeta_j^{(n)}\oplus 0)&=\zeta_j^{(n)}\oplus 0\quad \mbox{if $j\notin\Omega$}\,,\\
w_n(0\oplus \beta_j\ )  &=0\oplus\beta_j^{(n)}\quad \mbox{if $j\notin\Omega$}\,.
\end{align*}
Consider the new sequence $s_n=w_n\widehat{z}_n w_n^*$. Let $\puno$, $\pdos$, $\ptres$ and $\pcuatro$ be the orthogonal projections such that:
\begin{align*}
R(\puno)&=\mbox{span}\{\zeta_j^{(n)}\oplus 0:\, j\in\Omega\}\,,\\
R(\pdos)&=\mbox{span}\{\zeta_j^{(n)}\oplus 0:\, j\notin\Omega\}\,,\\
R(\ptres)&=\mbox{span}\{0\oplus\beta_j:\, j\in\Omega\}\,,\\
R(\pcuatro)&=\mbox{span}\{0\oplus\beta_j:\, j\notin\Omega\}\,.
\end{align*}
Note that, for each $n\in\N$, the operator $s_n$ commutes with the four projections.% $\puno$, $\pdos$,  $\ptres$, and $\pcuatro$.

\medskip

  \textbf{Claim:} there exists $n_0$ large enough so that
\begin{enumerate}
	\item $s_n (\puno+\pcuatro)= 0$ for every $n\in\N$;
	\item $\{\,|s_n\ptres|\,\}_{n\in\N}$ converges to an operator whose spectrum is contained in $\{2k\pi:\ k\in\Z\}$.
	\item There exists $\eps>0$ such that
	$$\Big(\bigcup_{n\ge n_0} \sigma(s_n \pdos) \Big)\cap \Big(\bigcup_{k\in\Z,k\neq 0} (2k\pi-\eps,2k\pi+\eps)\Big)=\varnothing\,.$$
	
\end{enumerate}

\medskip

  The first item is clear, and the second item is a direct consequence of \eqref{convergen a 2pi}. In order to prove the third one, recall that for every $j>M$ and every $n\in\N$ the moduli of the eigenvalues $\la_j(|z_n|)$ are contained in $(-\pi,\pi)$. On the other hand, we can take $n_0$ large enough so that the sequences $\{\la_j(z_n)\}_{n\in\N}$ for $j\in\{1,\ldots,M\}$ are close to their limits. Note that, for $j\notin\Omega$ the limits are far from the integer multiples of $2\pi$. These facts, all together, imply (3), and conclude the proof of the claim.

\medskip

  Straightforward computations show  that
$$
e^{i\widehat{z}}=\lim_{n\to\infty}e^{i(\widehat{u}_n\widehat{x}\widehat{u}_n^*+\widehat{v}_n\widehat{y}\widehat{v}_n^*)}=
\lim_{n\to\infty}e^{iw_n(\widehat{u}_n\widehat{x}\widehat{u}_n^*+\widehat{v}_n\widehat{y}\widehat{v}_n^*)w_n^*}\,,
$$
which implies
\begin{equation}\label{converge el truncadito}
e^{i\widehat{z}}=\lim_{n\to\infty}e^{i(s_n\pdos)}\,,
\end{equation}
because 
$$
\lim_{n\to\infty}e^{i(s_n(\ptres+\puno+\pcuatro))}=1. 
$$
The identity \eqref{converge el truncadito} and the claim allow us to apply Lemma \ref{rodriguez} to the sequence $\{s_n \pdos\}_{n\in\N}$, and to obtain a convergent subsequence. Therefore, the sequence $\{s_n\}_{n\in\N}$  has a convergent subsequence $\{s_{n_k}\}_{k\in\N}$. Let $s$ be its limit, that is
\begin{equation}\label{casi casi}
s=\lim_{k\to\infty} s_{n_k}=\lim_{k\to\infty} w_{n_k}(\widehat{u}_{n_k}\widehat{x}\,\widehat{u}_{n_k}^*+\widehat{v}_{n_k}\widehat{y}\,\widehat{v}_{n_k}^*)w_{n_k}^*\,.
\end{equation}
Clearly, this limit satisfies the identity $e^{i\widehat{z}}=e^{is}$.  On the other hand, if we consider the restriction of \eqref{casi casi} to $\ese=R(1-\pcuatro)$, then by Lemma \ref{decadentes} there are operators $\bar{s},\bar{x},\bar{y}\in B(\ese)$ which have the same non-zero eigenvalues (counted with multiplicity) as the operators $s$, $\widehat{x}$, and $\widehat{y}$. Extended as zero in $\ese^\bot$ (and using this notation),  $\bar{s}$, $\bar{x}$ and $\bar{y}$ become unitary equivalent to $s$, $\widehat{x}$, and $\widehat{y}$ respectively. Therefore, as $\bar{s}=\bar{x}+\bar{y}$, there exist two unitary operators $u_0$ and $v_0$ acting on $\cK$ such that
$$
s=u_0\widehat{x}u_0^*+v_0\widehat{y}v_0^*.
$$
This concludes the proof.
\edem

\section{Thompson-type formulae for operators in an embeddable II$_1$ factor}\label{fini}

  Throughout this section, let $\eme$ be a  II$_1$ factor that can be embedded in $\ultrar$.  We start with two technical lemmas.

\begin{lem}\label{aprox e}
Let $a,b\in\eme$ be Hermitian elements, and let $\{(A^{(m)},B^{(m)})\}_{m\in\N}$ be a sequence of matricial approximations. Then, for every polynomial $p\in\C[z,\bar{z}]$
$$
\tau \big( p(e^{ia}\cdot e^{ib})\big)=\lim_{m\to \infty} \tau_{n_m} \big(p(e^{iA^{(m)}}\cdot e^{iB^{(m)}})\big).
$$
\end{lem}
\bdem
It is a straightforward consequence of Theorem \ref{aprox}.
\edem

 Let us recall the definition of decreasing rearrangements of functions: given a measurable function $f:[0,1)\to\R$, its decreasing rearrangement $f^*:[0,1)\to\R$ is defined by
$$
f^*(t)=\inf\{s:\ |\{x:\,f(x)>s\}|\leq t\}\,.
$$
\begin{rem}\label{observacion reordenada}
Note that, given two functions $f,g:[0,1)\to\R$, if they satisfy
$|\{x:\,f(x)>s\}|=|\{x:\,g(x)>s\}|$ for every $s\in\R$, then $f^*=g^*$. The reader is referred to \cite{[BS]} for more details on decreasing rearrangements.
\EOE
\end{rem}

\begin{lem}\label{mentirita}
Let $f,g:[0,1)\to \R$ be bounded non-increasing functions such that $\|g\|_\infty\leq\pi$, and for any interval $I$ of the unit circle $S^1$ it holds that
\begin{equation}\label{trulala}
\int_0^1 \sub{\chi}{I}(e^{if(t)})\,dt=\int_0^1 \sub{\chi}{I}(e^{ig(t)})\,dt\,.
\end{equation}
Then, there is a function $\bar{g}:[0,1)\to\R$ such that $e^{if(t)}=e^{i\overline{g}(t)}$, and $\bar{g}^{\, *}=g$.
\end{lem}
\bdem
Let $\Omega=\{t\in[0,1):\, e^{if(t)}=-1\}$, and divide it in two measurable sets $\Omega_+$ and $\Omega_-$ such that 
$$
|\Omega_+|=|\{t\in[0,1):\ g(t)=\pi\}|\peso{and} |\Omega_-|=|\{t\in[0,1):\ g(t)=-\pi\}|\,.
$$
This is possible because $|\Omega|=|\{t\in[0,1): e^{ig(t)}=-1\}|$ by \eqref{trulala}.
Define $\bar{g}:[0,1)\to\R$ as follows:
$$ 
\bar{g}(t):=\begin{cases}
f(t)-2k\pi & \mbox{if $f(t)\in \big((2k-1)\pi, (2k+1)\pi\big)$;}\\
\ \ \,\pi& \mbox{if $t\in \Omega_+$;}\\
-\pi& \mbox{if $t\in \Omega_-$;}
\end{cases}\,.
$$
The function $\bar{g}$ clearly satisfies the identity $e^{if(t)}=e^{i\overline{g}(t)}$. So, for every arc $I$ of the unit circle
$$
\int_0^1 \sub{\chi}{I}(e^{i\overline{g}(t)})\,dt=\int_0^1 \sub{\chi}{I}(e^{if(t)})\,dt\,,
$$
and therefore
\begin{equation}\label{pares}
\int_0^1 \sub{\chi}{I}(e^{i\overline{g}(t)})\,dt=\int_0^1 \sub{\chi}{I}(e^{ig(t)})\,dt\,.
\end{equation}

% Note that, the a priori estimate $\|\bar{g}\|_\infty\leq \pi$ can be improved using \eqref{pares}. Indeed, straightforward computations show that $\|\bar{g}\|_\infty=\|g\|_\infty<\pi$.

\medskip

The next (and last) step, is to prove that $\bar{g}^*=g^*=g$ (almost everywhere). The last identity holds because $g$ is decreasing and the decreasing rearrangements considered here are with respect to the Lebesgue measure. To prove that $\bar{g}^*=g^*$, it is enough to verify that for every $s\in\R$
\begin{equation}\label{trulala2}
|\{x:\,\bar{g}(x)>s\}|=|\{x:\,g(x)>s\}|\,.
\end{equation}
Note that, by construction, $\|\bar{g}\|_\infty\leq \pi$. Hence, $\|\bar{g}\|_\infty=\|g\|_\infty$ by (\ref{pares}). Moreover, also by construction, it holds that
\begin{align*}
|\{x:\,\bar{g}(x)=-\pi\}|&=|\{x:\,g(x)=-\pi\}|.
\end{align*}
Therefore,  the equality in \eqref{trulala2}  is apparent if $s>\|g\|_\infty$ or $s\leq -\pi$.
On the other hand, if $-\pi< s\leq \|g\|_\infty$, let $I=\{e^{it}: s<t\leq\pi\}$. Then
\begin{align*}
|\{x:\,\bar{g}(x)>s\}|&=\int_0^1 \sub{\chi}{I}(e^{i\overline{g}(t)})\,dt - |\{x:\,\bar{g}(x)=-\pi\}|\\
&=\int_0^1 \sub{\chi}{I}(e^{ig(t)})\,dt - |\{x:\, g(x)=-\pi\}| \\
&=|\{x:\,{g}(x)>s\}|\,.
\end{align*}
This concludes the proof.
\edem

\begin{teo}\label{Thompson en el ultranoseque}
Given $a,b\in\eme$ Hermitian, there are two sequences of unitaries $\{u_n\}_{n\in\N}$ and $\{v_n\}_{n\in\N}$ such that
$$
e^{ia}\,e^{ib}=\lim_{n\to\infty} e^{i(u_nau_n^* + v_nbv_n^*)}\,,
$$
where the convergence is with respect to the operator norm topology.
\end{teo}

\bdem
Let $\{A^{(m)}\}_{m\in\N}$ and $\{B^{(m)}\}_{m\in\N}$ be sequences of matricial approximations of $a$ and $b$ respectively. By Thompson's theorem, there are unitary matrices $U_m$ y $V_m$ such that for each $m\in\N$
$$
e^{iA^{(m)}}\,e^{iB^{(m)}}=e^{i(U_mA^{(m)}U_m^*+V_mB^{(m)}B_m^*)}\,.
$$
Define $D_m=U_mA^{(m)}U_m^*+V_mB^{(m)}B_m^*$. By Theorem \ref{Horn finito}, the functions $\la_{A^{(m)}}$, $\la_{B^{(m)}}$, and $\la_{D^{(m)}}$  satisfy equation \eqref{integrohorn}. Since the sequence of non-increasing functions $\{\la_{D^{(m)}}\}_{m\in\N}$ is uniformly bounded, by Helly's selection theorem, there is a subsequence of this sequence that converges for all but almost countable many points $t\in[0,1)$. To simplify the notation, let us assume that the original sequence converges in this way, and let $f$ be its limit. This limit is also non-increasing and bounded. Moreover, as $\la_{A^{(m)}}$, $\la_{B^{(m)}}$, and $\la_{D^{(m)}}$ satisfy \eqref{integrohorn} for every $m\in\N$, by the dominated convergence theorem, $\la_a$, $\la_b$ and $f$ also satisfy those inequalities. Then, there are operators $a'$, $b'$ such that
\begin{equation}\label{los primos}
\la_{a'}=\la_{a}\ ,\quad \la_{b'}=\la_{b}\ ,\peso{and} \la_{a'+b'}=f
\end{equation}

Let $c\in\eme$ such that $e^{ia}\,e^{ib}=e^{ic}$ and $\|c\|\leq\pi$. Given a polynomial $p$, on one hand by Lemma \ref{aprox e}:
\begin{align}
\lim_{m\to\infty} \tau_{n_m} \big(p(e^{iA^{(m)}}\cdot e^{iB^{(m)}})\big)
&=\tau \big( p(e^{ia}\cdot e^{ib})\big)=\tau \big( p(e^{ic})\big)
=\int_0^1 p(e^{i\la_c(t)})\,dt\,. \label{con lambda c}\\
\intertext{ On the other hand, by the dominated convergence theorem, we obtain}
\lim_{m\to\infty} \tau_{n_m} \big(p(e^{iA^{(m)}}\cdot e^{iB^{(m)}})\big)&=\lim_{m\to\infty}\tau_{n_m} \big( p(e^{iD^{(m)}})\big)=
\lim_{m\to\infty}\int_0^1 p(e^{i\la_{D^{(m)}}(t)})\,dt \nonumber\\
&=\int_0^1 p(e^{if})\,dt\,. \label{con f}
\end{align}

 Therefore, \eqref{con lambda c} and \eqref{con f} imply that for every polynomial $p$
$$
\int_0^1 p(e^{i\la_c(t)})\,dt=\int_0^1 p(e^{if(t)})\,dt\,.
$$
Using standard arguments we obtain the same result replacing the polynomials by characteristic functions of arcs. Then, by Lemma \ref{mentirita}, there is a function $\bar{\la_c}$ such that $e^{if}=e^{i\overline{\la}_c}$, and $\bar{\la}_c^*=\la_c$. Suppose that
$$
c=\int_0^1 \la_c(t)\ de(t)\,,
$$
and define
\begin{align*}
c'=\int_0^1 \bar{\la}_c(t)\ de(t)\peso{and}
d=\int_0^1 f(t)\ de(t)\,.
\end{align*}
Then, $e^{ic'}=e^{id}$, $\la_d=f$, and $\la_{c'}=\la_c$. Combining these facts with equation \eqref{los primos}, and using Proposition \ref{Kamei}, we get that there are sequences  $\{u_n^{(a)}\}_{n\in\N}$, $\{u_n^{(b)}\}_{n\in\N}$, $\{u_n^{(c)}\}_{n\in\N}$, and $\{u_n^{(d)}\}_{n\in\N}$ of unitary elements of $\eme$ so that
\begin{align*}
d&=\lim_{n\to\infty} (u_n^{(d)})\,(a'+b')\,(u_n^{(d)})^*\ ,\\
a'&=\lim_{n\to\infty} (u_n^{(a)})\,a\,(u_n^{(a)})^*\ ,\\
b'&=\lim_{n\to\infty} (u_n^{(b)})\,b\,(u_n^{(b)})^*\ ,\ \ \mbox{and}\\
c&=\lim_{n\to\infty} (u_n^{(c)})\,c'\,(u_n^{(c)})^*\,.
\end{align*}
Finally, if we define $u_n=u_n^{(c)}u_n^{(d)}u_n^{(a)}$ y $v_n=u_n^{(c)}u_n^{(d)}u_n^{(b)}$ we get
$$
e^{ic}=\lim_{n\to\infty}e^{i(u_nau_n^* + v_nbv_n^*)}\,,
$$
which concludes the proof.
\edem

\appendix{

\section{Brief review on Horn's conjecture}\label{appendix A}

One of the most challenging problems in linear algebra has been to characterize the real $n$-tuples $\alpha$, $\beta$, and $\gamma$ that are the eigenvalues of $n\times n$ Hermitian matrices $A$, $B$, and $C$ such that $C=A+B$. In his remarkable 1962 paper \cite{Horn}, Roger Horn found necessary condition on the n-tuples $\alpha$, $\beta$, and $\gamma$ and conjectured that this conditions were also sufficient. This conjecture remained open for  several years, and it was solved at the end of the 20th century. Later on, these results were extended to operators in embeddable II$_1$ factors. In this appendix, we briefly recall these results; for some really deep material on the subject, we point the reader to the nice surveys by R. Bhatia \cite{[B]} and W. Fulton \cite{[F]}.

To begin with, we are going to fix some notation and conventions in order to state correctly the results in the finite dimensional setting. Firstly, the $n$-tuples will be considered arranged in non-increasing order, and by means of $\la(A)$ we denote the vector of eigenvalues of a self-adjoint matrix, also arranged in non-increasing order.

Clearly, one necessary condition that three $n$-tuples $\alpha$, $\beta$, and $\gamma$ have to satisfy in order to be the eigenvalues of Hermitian matrices $A$, $B$, and  $C$ such that $C=A+B$, is the next identity
\begin{equation}\label{igualito chico}
\sum_{j=1}^n \gamma_j=\sum_{j=1}^n \alpha_j + \sum_{j=1}^n \beta_j\,.
\end{equation}
This equality is far from being sufficient. In \cite{Horn}, Horn prescribed sets of triples $(I,J,K)$ of subsets of $\{1,\ldots,n\}$, that we will always write in increasing order, and he proved that the system of inequalities
$$
\sum_{k\in K}^n \gamma_k\leq \sum_{i\in I} \alpha_i + \sum_{j\in J} \beta_j\,,
$$
are necessary. The triples $(I,J,K)$ are defined by the following inductive procedure. Set
$$
U^n_r:=\left\{(I,J,K):\sum_{i\in I} i + \sum_{j\in J} j=\frac{r(r+1)}{2}+\sum_{k\in K} k \right\}\,.
$$
For $r=1$ set $T^n_1=U^n_1$. If $n\geq 2$, set
\begin{align*}
T^n_r:=\Big\{(I,J,K)\in U^n_r: \ &\mbox{for all $p<r$ and all $(F,G,H)\in T^r_p$, } \\
                             &\sum_{f\in F} i_f + \sum_{g\in G} j_g\leq \frac{p(p+1)}{2}+ \sum_{h\in H} k_h \Big\}\,.
\end{align*}
Then, the system of inequalities considered by Horn runs over all the triples in the set $\ete_n:=\bigcup_{k=1}^n T^n_k$. He also conjectured that this system of inequalities, together with the identity \eqref{igualito chico}, were sufficient. The proof of this conjecture is a consequence of several deep works of Klyachko, Knutson, and Tao (see \cite{[F],[Kl1],[KT],[KTW]}).

\begin{teo}\label{Horn finito}
Given $\alpha,\beta,\gamma\in\R^n$, the following statements are equivalent:
\begin{enumerate}
	\item There are $n\times n$ Hermitian matrices $A$, $B$ and $C$ such that $C=A+B$ and $\la(A)=\alpha$, $\la(B)=\beta$, and $\la(C)=\gamma$;
	\item $\sum_{k=1}^n\gamma_{k}= \sum_{i=1}^n\alpha_{i} + \sum_{j=1}^n \beta_{j}$, and for every $(I,J,K)$ in $T^n_r$, $\sum_{k\in K}\gamma_{k}\leq \sum_{i\in I}\alpha_{i} + \sum_{j\in J} \beta_{j}$
\end{enumerate}
\end{teo}

Later on, this result was extended by Bercovici and Li in \cite{[BL]} to  operators in an embeddable II$_1$ factor $\mathcal M$, i.e. a factor that can be embedded in the ultrapower of the hyperfinite factor. To state correctly this generalization, we need to introduce some notations.  Given $n\in\N$, if $I\subseteq \{1,\ldots,n\}$, then $\sigma_I$ denotes the set
$$
\bigcup_{i\in I}\left[\frac{(i-1)}{n},\frac{i}{n}\right)\,.
$$
With this notation, the set $\ete$ is defined by
$$
\displaystyle \ete:=\bigcup_{n=1}^\infty \bigcup_{r=1}^{n-1} \{(\sigma_I,\sigma_J,\sigma_K):\ (I,J,K)\in T^n_r\}\,.
$$
\begin{teo}\label{Horn II1}
Consider bounded non-increasing right-continuous functions $u$, $v$, and $w$ defined in the $[0,1)$. The following are equivalent:
\begin{enumerate}
	\item There exist $a,b\in\eme$ such that $u=\la_a$, $v=\la_b$ and $w=\la_{a+b}$;
	\item The functions $u$, $v$, and $w$ satisfy:
	\begin{align}
	&\int_0^1 u(t)\,dt+\int_0^1 v(t)\,dt=\int_0^1 w(t)\,dt\,,\ \  and \nonumber \\
	&\int_{\omega_1} u(t)\,dt+\int_{\omega_2} v(t)\,dt\geq \int_{\omega_3}
	 w(t)\,dt\,,\ \quad \forall (\omega_1,\omega_2,\omega_3)\in \ete\,. \label{integrohorn}
	 \end{align}
\end{enumerate}
\end{teo}
}

\bigskip

\fontsize {10}{10}\selectfont

Gabriel Larotonda and Alejandro Varela:\\
Instituto de Ciencias \\
Universidad Nacional de General Sarmiento. \\
J. M. Gutiérrez 1150 \\
(B1613GSX) Los Polvorines, \\
Buenos Aires, Argentina.  \\
e-mails: glaroton@ungs.edu.ar,\\
avarela@ungs.edu.ar\\

Jorge Antezana:\\
Universidad Nacional de La Plata.\\
Departamento de Matemática,\\
Esq. 50 y 115 s/n\\
Facultad de Ciencias\\
La Plata  (1900) \\
Buenos Aires, Argentina.\\
e-mail: antezana@mate.unlp.edu.ar\\

 J. Antezana, G. Larotonda\\
 and A. Varela:\\
 Instituto Argentino de Matemática\\
 ``Alberto Calder\'on'', CONICET\\
 Saavedra 15, 3er piso\\
 (C1083ACA) Buenos Aires,\\
 Argentina.\\

\end{document}